\input amstex
\documentstyle{amsppt}
\magnification=\magstep1
 \hsize 13cm \vsize 18.35cm \pageno=1
\loadbold \loadmsam
    \loadmsbm
    \UseAMSsymbols
\topmatter
\NoRunningHeads
\title Euler numbers and polynomials of higher order
\endtitle
\author
  T. Kim
\endauthor

\abstract The purpose of this paper is to present a systemic study
of some families of higher-order $q$-Euler numbers and polynomials
and we construct $q$-zeta function of order $r$ which interpolates
higher-order $q$-Euler numbers at negative integer.
\endabstract
\thanks  2000 AMS Subject Classification: 11B68, 11S80
\newline  The present Research has been conducted by the research
Grant of Kwangwoon University in 2010
\endthanks
\endtopmatter

\document

{\bf\centerline {\S 1. Introduction/ Preliminaries}}

 \vskip 15pt
Let $p$ be a fixed odd prime number. Throughout this paper $\Bbb
Z_p$, $\Bbb Q_p$, $\Bbb C$ and $\Bbb C_p$ will, respectively,
denote the ring of $p$-adic rational integers, the field of
$p$-adic rational numbers, the complex number field and the
completion of algebraic closure of $\Bbb Q_p$. Let $v_p$ be the
normalized exponential valuation of $\Bbb C_p$ with
$|p|_p=p^{-v_p(p)}=\frac{1}{p}$ (see [16]). When one talks of $q$-extension,
$q$ is variously considered  as an indeterminate, a complex number
$q\in \Bbb C$ or $p$-adic number $q\in \Bbb C_p$. If $q\in \Bbb
C$, one normally assumes $|q|<1$. If $q\in\Bbb C_p$, one normally
assumes $|1-q|_p<1.$ For a fixed $d\in \Bbb N$ with $(p,d)=1$, $
d\equiv 1$ $ ( mod \ 2)$ , we set
$$\aligned
&X =X_d =\varprojlim_N \Bbb Z/dp^N , X^\ast =\bigcup_{\Sb  0<a
<dp\\(a,p)=1
\endSb}a +d p\Bbb Z_p ,\\
& a +d p^N \Bbb Z_p = \{ x\in X | x \equiv a
\pmod{p^N}\},\endaligned$$ where $a\in\Bbb Z$ satisfies the
condition $0\leq a < d p^N$.
 The binomial formulae are known as
$$(1-b)^n=\sum_{i=0}^n \binom{n}{i} (-1)^ib^i, \text{ where $\binom
{n}{l}=\frac{ n (n-1)\cdots (n-l+1)}{l!}$},$$ and
$$\frac{1}{(1-b)^n}=(1-b)^{-n}=\sum_{l=0}^{\infty}\binom{-n}{l}(-b)^l
=\sum_{i=0}^{\infty}\binom{n+i-1}{i} b^i.$$ Recently, many authors
have studied the $q$-extension in the various area(see [4, 5, 6]).
In this paper, we try to consider the theory of $q$-integrals in
the $p$-adic number field associated with Euler numbers and
polynomials closely related to fermionic distribution. We say that
$f$ is uniformly differentiable function at a point $a\in\Bbb
Z_p$, and write $f\in UD(\Bbb Z_p),$ if the difference quotient
$F_{f}(x, y)=\frac{f(x)-f(y)}{x-y}$ have a limit $f^{\prime}(a)$
as $(x,y)\rightarrow (a,a).$ For $f \in UD(\Bbb Z_p),$ the
fermionic $p$-adic $q$-integral on $\Bbb Z_p$ is defined as
$$I_{q}(f)=\int_{\Bbb Z_p}f(x) d\mu_{q}(x)=\lim_{N\rightarrow
\infty}\frac{1+q}{1+q^{p^N}}\sum_{x=0}^{p^N-1}f(x)(-q)^x, \text{
(see [7, 8, 9, 16])}. \tag1$$ Thus, we note that
$$\lim_{q\rightarrow 1}I_{q}(f)=I_{1}(f)=\int_{\Bbb Z_p} f(x)
d\mu_{1}(x). \tag2$$ For $n\in \Bbb N$, let $f_n(x)= f(x+n)$. Then
we have
$$I_{1}(f_n)=(-1)^nI_{1}(f)+2\sum_{l=0}^{n-1}(-1)^{n-1-l}f(l). \tag3$$

Using formula (3), we can readily derive  the Euler polynomials,
$E_n(x),$ namely,
$$\int_{\Bbb Z_p} e^{(x+y)t}d\mu_{1}(y)=\frac{2}{e^t
+1}e^{xt}=\sum_{n=0}^{\infty}E_n(x)\frac{t^n}{n!}, \text{ (see
[16-20])}.$$ In the special case $x=0$, the sequence
$E_n(0)=E_n$ are called the $n$-th Euler numbers. In one of an
impressive series of papers( see[1, 2, 3, 21, 23]), Barnes developed the
so-called multiple zeta and multiple gamma functions. Barnes'
multiple zeta function $\zeta_N (s, w| a_1, \cdots, a_N)$ depend
on the parameters $a_1, \cdots, a_N$ that will be assumed to be
positive. It is defined by the following series:
$$\zeta_{N}(s, w|a_1, \cdots, a_N)=\sum_{m_1, \cdots,
m_N=0}^{\infty}(w+m_1a_1+\cdots+m_Na_N)^{-s}\text{ for $\Re(s)>N,
\Re(w)>0 .$}\tag4$$ From (4), we can easily see that
$$\zeta_{M+1}(s, w+a_{M+1}|a_1, \cdots, a_{N+1})-\zeta_{M+1}(s, w|a_1,
\cdots, a_{N+1})=-\zeta_{M}(s, w|a_1, \cdots, a_N),$$
 and $\zeta_0(s, w)=w^{-s}$( see [1]). Barnes showed that
 $\zeta_N$ has a meromorphic continuation in $s$ (with simple
 poles only at $s=1, 2, \cdots, N$ and defined his multiple gamma
 function $\Gamma_N(w)$ in terms of the $s$-derivative at $s=0,$
 which may be recalled here as follows:   $\psi_n(w|a_1,\cdots, a_N)=
\partial_s\zeta_N(s,w|a_1, \cdots, a_N)|_{s=0}$( see[11]).
Barnes' multiple Bernoulli polynomials $B_n(x,r|a_1, \cdots, a_r)$
are defined by
$$\frac{t^r}{\prod_{j=1}^r(e^{a_jt}-1)}e^{xt}=\sum_{n=0}^{\infty}B_n(x,
r|a_1, \cdots, a_r)\frac{t^n}{n!}, \text{ ($|t|<\max_{1\leq i\leq r}\frac{2\pi}{|a_i|}$)}, \text{ (see [1, 11])}.\tag5$$ By
(4) and (5), we see that
$$\zeta_N(-m, w|a_1, \cdots, a_N)=\frac{(-1)^N
m!}{(N+m)!}B_{N+m}(w,N|a_1, \cdots, a_N), \text{ (see [1])}, $$ where $w>0$ and $m$
is a positive integer. By using the fermionic $p$-adic
$q$-integral on $\Bbb Z_p$, we consider the Barnes' type multiple
$q$-Euler polynomials and numbers in this paper. The main purpose
of this paper is to study a systemic properties of some families
of higher-order $q$-Euler polynomials and numbers. Finally, we
construct $q$-zeta function of order $r$ which interpolates
higher-order $q$-Euler numbers and polynomials at negative
integer.

\vskip 10pt

{\bf\centerline {\S 2. higher-order $q$-Euler numbers and
polynomials }} \vskip 10pt

Let $x, w_1, w_2, \cdots, w_r$ be complex numbers with positive
real parts. In $\Bbb C$, the Barnes type multiple Euler numbers
and polynomials are defined by
$$\frac{2^r}{\prod_{j=1}^r(e^{w_jt}+1)}e^{xt}=\sum_{n=0}^{\infty}E_{n}^{(r)}(x|w_1,
\cdots, w_r)\frac{t^n}{n!}, \text{ for $|t|< max \{
\frac{\pi}{|w_i|}|i=1,\cdots,r \}$},\tag6$$ and $E_n^{(r)}(w_1,
\cdots, w_r)=E_n^{(r)}(0|w_1, \cdots, w_r)$(see [11, 12, 14]). In this
section, we assume that $q\in \Bbb C_p$ with $|1-q|_p<1$. We first
consider the $q$-extension of Euler polynomials as follows:
$$\sum_{n=0}^{\infty}E_{n,q}(x)\frac{t^n}{n!}=\int_{\Bbb
Z_p}q^ye^{(x+y)t}d\mu_{1}(y)=2\sum_{m=0}^{\infty}(-q)^me^{(m+x)t}=\frac{2}{qe^t+1}e^{xt}
.\tag7$$ In the special case $x=0$, $E_{n,q}=E_{n,q}(0)$ are
called the $q$-Euler numbers. By using multivariate $p$-adic
invariant integral on $\Bbb Z_p$, we consider the $q$-Euler
polynomials of order $r\in\Bbb N$ as follows:
$$\aligned
&\sum_{n=0}^{\infty}E_{n,q}^{(r)}(x)\frac{t^n}{n!}=\int_{\Bbb
Z_p}\cdots\int_{\Bbb
Z_p}e^{(x+x_1+\cdots+x_r)t}q^{x_1+\cdots+x_r}d\mu_{1}(x_1)\cdots d\mu_{1}(x_r)\\
&=\left(\frac{2}{qe^t+1}\right)^r e^{xt}
=2^r\sum_{m=0}^{\infty}\binom{m+r-1}{m}(-q)^m e^{(m+x)t}.
\endaligned\tag8$$
In the special case $x=0$, the sequence
$E_{n,q}^{(r)}(0)=E_{n,q}^{(r)}$ are refereed as the $q$-extension
of the Euler numbers of order $r$. Let $f\in\Bbb N$ with $f\equiv
1$ $(mod  \ 2)$. Then we have
$$\aligned
&E_{n,q}^{(r)}(x)=\int_{\Bbb Z_p}\cdots\int_{\Bbb
Z_p}q^{x_1+\cdots+x_r}(x+x_1+\cdots+x_r)^n d\mu_{1}(x_1)\cdots d\mu_{1}(x_r)\\
 &=2^r\sum_{m_1, \cdots,
m_r=0}^{\infty}(-q)^{m_1+\cdots+m_r}(m_1+\cdots+m_r+x)^n.
\endaligned\tag9$$
By (8) and (9), we obtain the following theorem.

\proclaim{ Theorem 1} For $n\in \Bbb Z_+$, we have
$$\aligned
 E_{n,q}^{(r)}(x)&=2^r\sum_{m_1, \cdots,
m_r=0}^{\infty}(-q)^{m_1+\cdots+m_r}(m_1+\cdots+m_r+x)^n\\
&=2^r\sum_{m=0}^{\infty}\binom{m+r-1}{m}(-q)^m (m+x)^n.
\endaligned$$
\endproclaim
Let
$F_q^{(r)}(t,x)=\sum_{n=0}^{\infty}E_{n,q}^{(r)}(x)\frac{t^n}{n!}.$
Then we have
$$\aligned
F_q^{(r)}(t,x)&=2^r\sum_{m=0}^{\infty}\binom{m+r-1}{m}(-q)^me^{(m+x)t}\\
&=2^r\sum_{m_1, \cdots,
m_r=0}^{\infty}(-q)^{m_1+\cdots+m_r}e^{(m_1+\cdots+m_r+x)t}.
\endaligned\tag10$$
Let $\chi$ be the Dirichlet's character with conductor $f\in\Bbb
N$ with $f\equiv 1$ $( mod   \ 2)$. Then the generalized $q$-Euler
polynomials attached to $\chi$ are defined by
$$\sum_{n=0}^{\infty}E_{n,\chi,q}(x)\frac{t^n}{n!}=2\sum_{m=0}^{\infty}(-q)^m\chi(m)e^{(m+x)t}.
\tag11$$

Thus, we have
$$\aligned
&E_{n,\chi,q}(x)\\
&=\sum_{a=0}^{f-1}\chi(a)(-q)^a\int_{\Bbb Z_p}(x+a+fy)^n q^{fy}
d\mu_{1}(y)=f^n\sum_{a=0}^{f-1}\chi(a)(-q)^aE_{n,
q^f}(\frac{x+a}{f}).\endaligned\tag12$$ In the special case $x=0$,
the sequence $E_{n,\chi,q}(0)=E_{n,\chi,q}$ are called the $n$-th
generalized $q$-Euler numbers attached to $\chi$. From (2) and
(3), we can easily derive the following equation.
$$E_{m,\chi,q}(nf)-(-1)^nE_{m,\chi,q}=2\sum_{l=0}^{nf-1}(-1)^{n-1-l}
\chi(l)q^l l^m.$$

Let us define higher-order generalized $q$-Euler polynomials
attached to $\chi$ as follows:
$$\aligned
&\int_X\cdots\int_X\left( \prod_{i=1}^r
\chi(x_i)\right)e^{(x_1+\cdots+x_r+x)t}q^{x_1+\cdots
+x_r}d\mu_{1}(x_1)\cdots d\mu_{1}(x_r)\\
&=\left(\frac{2\sum_{a=0}^{f-1}(-q)^a\chi(a)e^{at}}{q^fe^{ft}+1}
\right) =\sum_{n=0}^{\infty}E_{n,\chi,q}^{(r)}(x)\frac{t^n}{n!},
\endaligned\tag13$$ where $E_{n,\chi,q}^{(r)}(x)$ are called the $n$-th
generalized $q$-Euler polynomials of order $r$ attached to $\chi$.
By (13), we see that
$$\aligned
&E_{n,\chi,q}^{(r)}(x)\\
&=2^r\sum_{m=0}^{\infty}\binom{m+r-1}{m}(-q^f)^m\sum_{a_1,\cdots,
a_r=0}^{f-1}(\prod_{j=1}^r\chi(a_j))(-q)^{\sum_{i=1}^r
a_i}(\sum_{j=1}^r a_j+x+mf)^n,
\endaligned\tag14$$
and
$$\sum_{n=0}^{\infty}E_{n,\chi,q}^{(r)}(x)\frac{t^n}{n!}
=2^r\sum_{m_1,\cdots,
m_r=0}^{\infty}(-q)^{\sum_{j=1}^rm_j}\left(\prod_{i=1}^r\chi(m_i)\right)
e^{(x+\sum_{j=1}^rm_j)t}. \tag15$$ In the special case $x=0$, the
sequence $E_{n,\chi,q}^{(r)}(0)=E_{n,\chi,q}^{(r)}$ are called the
$n$-th generalized $q$-Euler numbers of order $r$ attached to
$\chi$.

By (14) and (15), we obtain the following theorem.

\proclaim{ Theorem 2} Let $\chi$ be the Dirichlet's character with
conductor $f\in\Bbb N$ with $f\equiv 1$ $(mod  \ 2)$.
 For $n\in \Bbb Z_+$, $r\in\Bbb N$,  we have
$$\aligned
& E_{n,\chi,q}^{(r)}(x)\\
&=2^r\sum_{m=0}^{\infty}\binom{m+r-1}{m}(-q^f)^m\sum_{a_1,\cdots,
a_r=0}^{f-1}(\prod_{j=1}^r\chi(a_j))(-q)^{\sum_{i=1}^r
a_i}(\sum_{j=1}^r a_j+x+mf)^n\\
&=2^r\sum_{m_1,\cdots,
m_r=0}^{\infty}(-q)^{m_1+\cdots+m_r}\left(\prod_{i=1}^r\chi(m_i)\right)
(x+m_1+\cdots+m_r)^n.
\endaligned$$
\endproclaim
  For $h\in\Bbb Z$ and $r\in\Bbb N$, we introduce the extended higher-order $q$-Euler polynomials
  as follows:
  $$E_{n,q}^{(h,r)}(x)=\int_{\Bbb Z_p}\cdots \int_{\Bbb Z_p}q^{\sum_{j=1}^r(h-j)x_j}(x+x_1+\cdots+x_r)^n
  d\mu_{1}(x_1)\cdots d\mu_1(x_r).\tag16$$
  From (16), we note that
  $$E_{n,q}^{(h,r)}(x)=2^r\sum_{m_1,\cdots, m_r=0}^{\infty}q^{(h-1)m_1+\cdots+(h-r)m_r}(-1)^{m_1+\cdots+m_r}
  (x+m_1+\cdots+m_r)^n,\tag17$$
  where ${\binom{n}{l}}_q=\frac{[n]_q [n-1]_q\cdots
  [n-l+1]_q}{[l]_q[l-1]_q\cdots[2]_q[1]_q}$ and
  $[n]_q=\frac{1-q^n}{1-q}.$

 Thus, we have
  $$ E_{n,q}^{(h,r)}(x)=2^r\sum_{m=0}^{\infty}{\binom{m+r-1}{m}}_q (-q^{h-r})^m(x+m)^n.\tag18$$
  Let
  $$\aligned
  F_q^{(h,r)}(t,x)&=\int_{\Bbb Z_p}\cdots \int_{\Bbb Z_p}
  q^{\sum_{j=1}^r(h-j)x_j}e^{(\sum_{i=1}^r x_i+x)t}d\mu_1(x_1)\cdots
  d\mu_1(x_r)\\
  &  =\sum_{n=0}^{\infty}E_{n,q}^{(h,r)}(x)\frac{t^n}{n!}.\endaligned$$
  Then we have
    $$\aligned
   F_q^{(h,r)}(t,x)&=\frac{2^r}{\prod_{j=1}^{r}(1+e^tq^{h-r+j-1})}e^{xt}
   =2^r\sum_{m=0}^{\infty}{\binom{m+r-1}{m}}_q(-q^{h-r})^me^{(m+x)t}\\
   &=2^r \sum_{m_1,\cdots, m_r=0}^{\infty}q^{\sum_{j=1}^r (h-j)m_j}(-1)^{\sum_{j=1}^r m_j}e^{(x+m_1+\cdots+m_r)t}.
    \endaligned\tag19$$
 Therefore, we obtain the following theorem.
  \proclaim{ Theorem 3} For $h, \in \Bbb Z$, $r\in\Bbb N$, and $x\in \Bbb Q^{+}$,
   we have
$$\aligned
 E_{n,q}^{(h,r)}(x)&=2^r\sum_{m_1, \cdots,
m_r=0}^{\infty}q^{(h-1)m_1+\cdots+(h-r)m_r}(-1)^{m_1+\cdots+m_r}(m_1+\cdots+m_r+x)^n\\
&=2^r\sum_{m=0}^{\infty}{\binom{m+r-1}{m}}_q (-q^{h-r})^m(x+m)^n.
\endaligned$$
\endproclaim
For $f\in\Bbb N$ with $f\equiv 1$ $(mod   \ 2)$, it is easy to
show that  the following distribution relation for
$E_{n,q}^{(h,r)}(x)$.
 $$E_{n,q}^{(h,r)}(x)=f^n\sum_{a_1,\cdots, a_r=0}^{f-1}(-1)^{a_1+\cdots+a_r}q^{\sum_{j=1}^r(h-j)a_j}
  E_{n,q^f}(\frac{x+a_1+\cdots+a_r}{f}).$$
Let us consider Barnes' type higher-order $q$-Euler polynomials.
For $w_1, \cdots, w_r \in\Bbb Z_p$, we define  the Barnes' type
$q$- Euler polynomials of order $r$ as follow:
$$\aligned
&\sum_{n=0}^{\infty}E_{n,q}^{(r)}(x|w_1, \cdots, w_r)
\frac{t^n}{n!}=\frac{2^r}{\prod_{i=1}^r(e^{w_it}q^{w_i}+1)}e^{xt}\\
&= \int_{\Bbb Z_p}\cdots \int_{\Bbb Z_p}e^{(\sum_{j=1}^r
w_jx_j+x)t}q^{w_1x_1+\cdots+w_rx_r} d\mu_{1}(x_1)\cdots
d\mu_{1}(x_r).\endaligned\tag20$$ From (20), we can easily derive
the following equation.
$$ E_{n,q}^{(r)}(x|w_1, \cdots, w_r)=
\int_{\Bbb Z_p}\cdots\int_{\Bbb Z_p}(\sum_{i=1}^r
x_iw_i+x)^nq^{w_1x_1+\cdots+x_rw_r}d\mu_1(x_1)\cdots d\mu_1(x_r) .
\tag21$$ Thus, we have
 $$\aligned
 &E_{n,q}^{(r)}(x|w_1, \cdots, w_r)\\
 &=f^n\sum_{a_1, \cdots, a_r=0}^{f-1}(-1)^{\sum_{i=1}^ra_i}
 q^{\sum_{j=1}^r w_ja_j}E_{n, q^f}^{(r)}(\frac{\sum_{j=1}^rw_ja_j+x}{f}|w_1, \cdots, w_r),
 \endaligned\tag22 $$
   where $f\in\Bbb N$ with $f\equiv 1   \  (mod    \   2)$.
   By (22), we see that
   $$ E_{n,q}^{(r)}(x|w_1, \cdots, w_r)=2^r\sum_{m_1,\cdots, m_r=0}^{\infty}(-q)^{m_1w_1+\cdots+m_rw_r}
   (x+w_1m_1+\cdots+w_rm_r)^n.\tag23$$
In the special case $x=0$, the sequence $E_{n,q}^{(r)}(w_1,
\cdots, w_r)=E_{n,q}^{(r)}(0|w_1, \cdots, w_r)$ are called the
$n$-th Barnes' type  $q$-Euler numbers of order $r$.

Let $F_q^{(r)}(t,x|w_1, \cdots,
w_r)=\sum_{n=0}^{\infty}E_{n,q}^{(r)}(x|w_1, \cdots,
w_r)\frac{t^n}{n!}.$ Then we have
$$F_q^{(r)}(t,x|w_1, \cdots, w_r)=2^r\sum_{m_1, \cdots, m_r=0}^{\infty}(-q)^{m_1w_1+\cdots+m_rw_r}
e^{(x+w_1m_1+\cdots+w_rm_r)t}. \tag24$$ Therefore we obtain the
following theorem.
   \proclaim{ Theorem 4} For $w_1, \cdots, w_r \in \Bbb Z_p$, $r\in\Bbb N$, and $x\in \Bbb Q^{+}$,
     we have
  $$  E_{n,q}^{(r)}(x| w_1, \cdots, w_r)=2^r\sum_{m_1, \cdots,
  m_r=0}^{\infty}(-q)^{m_1w_1+\cdots+m_rw_r}(x+m_1w_1+\cdots+m_rw_r)^n .$$
  \endproclaim

  For $w_1, \cdots, w_r \in \Bbb Z_p$, $a_1, \cdots, a_r \in \Bbb Z$, we consider another $q$-extension
  of Barnes' type $q$-Euler polynomials of order $r$ as follows:
 $$
 \sum_{n=0}^r E_{n,q}^{(r)}(x|w_1,\cdots, w_r;a_1,\cdots,
 a_r)\frac{t^n}{n!}
 =\frac{2^r}{(q^{a_1}e^{w_1t}+1)(q^{a_2}e^{w_2t}+1)\cdots
(q^{a_r}e^{w_rt}+1)}. \tag25$$
 Thus, we have
   $$\aligned
   &E_{n,q}^{(r)}(x|w_1,\cdots,w_r;a_1,\cdots, a_r)\\
   &=\int_{\Bbb Z_p}\cdots \int_{\Bbb Z_p}e^{(x+\sum_{j=1}^r w_jx_j)t}q^{\sum_{j=1}^ra_jx_j}
   d\mu_1(x_1)\cdots d\mu_1(x_r).\endaligned\tag26$$
 From (25) and (26), we note that
 $$E_{n,q}^{(r)}(x|w_1,\cdots, w_r;a_1, \cdots, a_r)
 =2^r\sum_{m_1,\cdots, m_r=0}(-1)^{\sum_{j=1}^rm_j}q^{\sum_{i=1}^r a_im_i}(x+\sum_{j=1}^rw_jx_j)^n .\tag27$$
 Let $F_q^{(r)}(t,x|w_1,\cdots, w_r;a_1,\cdots, a_r)
 =\sum_{n=0}^{\infty}E_{n,q}^{(r)}(x|w_1,\cdots, w_r;a_1, \cdots, a_r)\frac{t^n}{n!}.$
  Then, we see that
 $$\aligned
 &F_q^{(r)}(t,x|w_1,\cdots, w_r;a_1,\cdots, a_r)\\
 &=2^r\sum_{m_1,\cdots, m_r=0}^{\infty}
 (-1)^{m_1+\cdots+m_r}q^{ a_1m_1+\cdots+a_rm_r}e^{(x+w_1m_1 +\cdots+ w_rm_r
 )t}.\endaligned\tag28$$
   \proclaim{ Theorem 5} For  $r\in\Bbb N$, $w_1, \cdots, w_r\in\Bbb Z_p$, and $a_1, \cdots, a_r \in \Bbb Z$,
   we have
 $$ E_{n,q}^{(r)}(x|w_1, \cdots, w_r; a_1, \cdots a_r)=2^r\sum_{m_1, \cdots,
 m_r=0}^{\infty}(-1)^{\sum_{j=1}^rm_j}q^{\sum_{i=1}^r a_im_i}(x+\sum_{j=1}^r w_jm_j)^n.$$
 \endproclaim
 Let $\chi$ be a Dirichlet's character with conductor $f\in \Bbb N$ with $f\equiv 1$ $(mod  \  2)$.
By using multivariate $p$-adic invariant integral on $X$, we now
consider the generalized Barnes' type $q$-Euler polynomials of
order $r$ attached to $\chi$ as follows:
 $$\aligned
 &\sum_{n=0}^{\infty}E_{n,\chi,q}^{(r)}(x|w_1, \cdots, w_r;a_1, \cdots, a_r)\frac{t^n}{n!}\\
& =\int_X\cdots \int_X
e^{(x+w_1x_1+\cdots+w_rx_r)t}\left(\prod_{j=1}^r\chi(x_j)\right)q^{a_1x_1+\cdots
+a_rx_r}d\mu_{1}(x_1)\cdots d\mu_{1}(x_r) .\endaligned$$

 Thus, we have
$$\aligned
&\left(\frac{2\sum_{b_1=0}^{f-1}\chi(b_1)q^{a_1b_1}(-1)^{b_1}e^{w_1b_1t}}{q^{a_1f}e^{w_1ft}+1}\right)
\times\cdots \times
\left(\frac{\sum_{b_r=0}^{f-1}\chi(b_r)q^{a_rb_r}(-1)^{b_r}e^{w_rb_rt}}{q^{a_rf}e^{w_rft}+1}\right)\\
&=\sum_{n=0}^{\infty}E_{n, \chi,q}^{(r)}(x|w_1,\cdots, w_r;a_1,
\cdots, a_r)\frac{t^n}{n!}.  \endaligned\tag 29$$

   From (29), we have
    $$\aligned
    &E_{n, \chi,q}^{(r)}(x|w_1,\cdots, w_r;a_1, \cdots, a_r) \\
    &=2^r\sum_{m_1, \cdots, m_r=0}^{\infty}\left(\prod_{j=1}^r\chi(m_i)\right)(-1)^{m_1+\cdots+m_r}
    q^{a_1m_1+\cdots+a_rm_r}(x+\sum_{j=1}^rw_jm_j)^n    . \endaligned$$

Therefore we obtain the following theorem.
  \proclaim{ Theorem 6} For  $r\in\Bbb N$, $w_1, \cdots, w_r\in\Bbb Z_p$, and $a_1, \cdots, a_r \in \Bbb Z$,
   we have
 $$\aligned
    &E_{n, \chi,q}^{(r)}(x|w_1,\cdots, w_r;a_1, \cdots, a_r) \\
    &=2^r\sum_{m_1, \cdots, m_r=0}^{\infty}\left(\prod_{j=1}^r\chi(m_i)\right)(-1)^{m_1+\cdots+m_r}
    q^{a_1m_1+\cdots+a_rm_r}(x+\sum_{j=1}^rw_jm_j)^n    . \endaligned$$
 \endproclaim
 Let $F_{q, \chi}^{(r)}(t,x|w_1,\cdots, w_r;a_1,\cdots, a_r)
  =\sum_{n=0}^{\infty}E_{n,\chi, q}^{(r)}(x|w_1,\cdots, w_r;a_1, \cdots, a_r)\frac{t^n}{n!}.$

  By Theorem 6, we see that
  $$\aligned
   &F_{q, \chi}^{(r)}(t,x|w_1,\cdots, w_r;a_1,\cdots, a_r) \\
   &= 2^r\sum_{m_1, \cdots, m_r=0}^{\infty}\left(\prod_{j=1}^r\chi(m_i)\right)(-1)^{m_1+\cdots+m_r}
      q^{a_1m_1+\cdots+a_rm_r}e^{(x+\sum_{j=1}^rw_jm_j)t}.    \endaligned\tag30$$

\vskip 10pt
   {\bf\centerline {\S 3. Higher-order $q$-zeta functions in $\Bbb C$}} \vskip 10pt

   In this section, we assume that $q\in\Bbb C$ with $|q|<1$ and the parameters $w_1,
  \cdots, w_r$ are positive.  From (28), we can define the Barnes' type $q$-Euler polynomials of order $r$
  in $\Bbb C$ as follows:
  $$\aligned
   &F_q^{(r)}(t, x|w_1, \cdots, w_r)=\frac{2^r}{\prod_{j=1}^r(e^{w_jt}q^{w_j}+1)}\\
   &= 2^r\sum_{m_1,\cdots, m_r=0}^{\infty}(-1)^{m_1+\cdots+m_r}q^{w_1m_1+\cdots+w_rm_r}
   e^{(x+w_1m_1+\cdots+w_rm_r)t}\\
   &=\sum_{n=0}^{\infty} E_{n,q}^{(r)}(x|w_1, \cdots, w_r)\frac{t^n}{n!},
   \text{ for $|t+\ln q|<\max_{1\leq i \leq r} \{ \frac{\pi}{|w_i|}\}$ }.
 \endaligned\tag31$$
     For $s, x \in \Bbb C$ with $\Re(x)>0,$
  we can derive the following Eq.(32) from the  Mellin transformation of $F_q^{(r)}(t,x|w_1, \cdots, w_r)$.
  $$\aligned
  &\frac{1}{\Gamma(s)}\int_{0}^{\infty}t^{s-1}F_q^{(r)}(-t,x|w_1, \cdots, w_r)dt\\
 & =2^r \sum_{m_1, \cdots, m_r=0}^{\infty}\frac{(-1)^{m_1+\cdots+m_r}
 q^{ m_1w_1+\cdots+m_rw_r}}{(x+ w_1m_1+\cdots+w_rm_r)^s}.
 \endaligned \tag32$$
  For $s, x \in \Bbb C$ with $\Re(x)>0$,   we define Barnes' type $q$-zeta function of order $r$ as follows:
  $$\zeta_{q}^{(r)}(s,x|w_1, \cdots, w_r)
  =2^r\sum_{m_1, \cdots, m_r=0}^{\infty}\frac{(-1)^{m_1+\cdots+m_r}q^{m_1w_1+\cdots+m_rw_r}}{(x+w_1m_1+\cdots+w_rm_r)^s}
 .\tag33$$
 Note that
 $\zeta_{q}^{(r)}(s,x|w_1,\cdots, w_r)$  is meromorphic function in whole complex $s$-plane. By using the Mellin transformation and the Cauchy
 residue theorem, we obtain the following theorem.
  \proclaim{ Theorem 7} For $x \in \Bbb C$ with $\Re(x)>0$, $n\in\Bbb Z_{+}$, we have
   $$\zeta_{q}^{(r)}(-n, x|w_1, \cdots, w_r)
   =E_{n,q}^{(r)}(x|w_1, \cdots, w_r).$$
    \endproclaim
    Let $\chi$ be a Dirichlet's character with conductor $f\in \Bbb N$ with $f\equiv 1$ $(mod  \  2)$.
 From (30), we can define the generalized Barnes' type $q$-Euler polynomials of order $r$ attached
 to $\chi$ in $\Bbb C$ as follows:
 $$\aligned
   &F_{q, \chi}^{(r)}(t,x|w_1,\cdots, w_r)\\
   &=\left(\frac{2\sum_{b_1=0}^{f-1}\chi(b_1)q^{w_1b_1}(-1)^{b_1}e^{w_1b_1t}}{q^{w_1f}e^{w_1ft}+1}\right)
   \times\cdots \times
\left(\frac{\sum_{b_r=0}^{f-1}\chi(b_r)q^{w_rb_r}(-1)^{b_r}e^{w_rb_rt}}{q^{w_rf}e^{w_rft}+1}\right) \\
   &= 2^r\sum_{m_1, \cdots, m_r=0}^{\infty}\left(\prod_{j=1}^r\chi(m_i)\right)(-1)^{m_1+\cdots+m_r}
      q^{w_1m_1+\cdots+w_rm_r}e^{(x+\sum_{j=1}^rw_jm_j)t}\\
      &= \sum_{n=0}^{\infty} E_{n,\chi,q}^{(r)}(x|w_1, \cdots, w_r)\frac{t^n}{n!}.    \endaligned\tag34$$

       From (34) and  Mellin transformation of $F_{q, \chi}^{(r)}(t,x|w_1, \cdots, w_r)$,
        we can easily derive the following equation (35) .
  $$\aligned
  &\frac{1}{\Gamma(s)}\int_{0}^{\infty}t^{s-1}F_{q, \chi}^{(r)}(-t,x|w_1, \cdots, w_r)dt\\
 & = 2^r\sum_{m_1, \cdots, m_r=0}^{\infty}\frac{\left( \prod_{j=1}^r \chi(m_i)\right)(-1)^{m_1+\cdots+m_r}
 q^{ m_1w_1+\cdots+m_rw_r}}{(x+ w_1m_1+\cdots+w_rm_r)^s}.
 \endaligned \tag35$$

  For $s, x \in \Bbb C$ with $\Re(x)>0$,  we also define Dirichlet's type Euler $q$-$l$-function of order $r$
  as follows:
  $$\aligned
  &l_{q}^{(r)}(s, x;\chi|w_1, \cdots, w_r)\\
 & =2^r\sum_{m_1, \cdots, m_r=0}^{\infty}\frac{\left(\prod_{j=1}^r\chi(m_j)\right)(-1)^{m_1+\cdots+m_r}
 q^{m_1w_1+\cdots+m_rw_r}}{(x+w_1m_1+\cdots+w_rm_r)^s}
 .\endaligned\tag36$$
    Note that
 $l_{q}^{(r)}(s,x;\chi|w_1,\cdots, w_r)$  is meromorphic function in whole complex $s$-plane.
 By using (34), (35), (36), and the Cauchy
 residue theorem, we obtain the following theorem.

  \proclaim{ Theorem 8} For $x, s \in \Bbb C$ with $\Re(x)>0$, $n\in\Bbb Z_{+}$, we have
   $$l_{q}^{(r)}(-n, x;\chi|w_1, \cdots, w_r)
   =E_{n,\chi, q}^{(r)}(x|w_1, \cdots, w_r).$$
    \endproclaim
 We note that Theorem 8 is $r$-ple Dirichlet's type $q$-$l$-series.  Theorem 8 seems to be  interesting and worthwhile
 for doing study in the area of  multiple  $l$-function
 related to the number theory.

\vskip 20pt

\quad Taekyun Kim

\quad Division of General Education-Mathematics,

\quad Kwangwoon University,

\quad Seoul 139-701, S. Korea \quad e-mail:\text{
tkkim$\@$kw.ac.kr}

\enddocument